\documentclass[10pt]{amsart}

\usepackage{amsmath, amssymb, amsfonts}
\usepackage[mathscr]{euscript}

\keywords{greedy decomposition, normal decomposition, symmetric normal decomposition, $\Delta$-normal decomposition, Garside family, Garside map, Garside element, Garside monoid, Garside group, Word Problem, Conjugacy Problem, braid group, Artin--Tits group, Deligne-Luzstig variety, self-distributivity, ordered group, Yang--Baxter equation, cell decomposition}

\def\Part#1{\bigskip\noindent{\bf PART #1}\quad\dotfill\quad}
\def\Chapter#1{\goodbreak\bigskip\noindent {\bf #1}\quad\dotfill\quad}
\def\Section#1{\smallskip {\rm #1}\quad\dotfill\quad}
\def\Subsection#1{\indent\indent #1\quad\dotfill\quad}
\def\Subsections#1{\indent\indent #1}
\def\Exercises{\smallskip Exercises\quad\dotfill\quad}
\def\Notes{Notes\quad\dotfill\quad}

\newcommand\BP[1]{B_{#1}^{\scriptscriptstyle+}}
\newcommand\CCC{\mathcal{C}}
\newcommand\CCCu{\CCC_1}
\newcommand\HS[1]{\hspace{#1ex}}
\newcommand\LDM{M_{\LDsmall}}
\newcommand\LDsmall{{\HS{-0.3}\scriptscriptstyle L\!D}}
\newcommand\NNNN{\mathbb{N}}
\newcommand\pdots{\hspace{0.2ex}{\cdot}{\cdot}{\cdot}\hspace{0.2ex}}
\newcommand\sig[1]{\sigma_{\hspace{-0.2ex}#1}^{\null}}
\newcommand\SSS{\mathcal{S}}
\def\VR(#1,#2){\vrule width0pt height#1mm depth#2mm}
\newcommand\wdots{, ...\hspace{0.2ex},}
\newcommand\ZZZZ{\mathbb{Z}}

\begin{document}

\begin{center}
{\bf FOUNDATIONS OF GARSIDE THEORY} 

\vspace{5mm}

\textsc{Patrick DEHORNOY}

\medskip with

\textsc{Fran\c cois DIGNE}

\textsc{Eddy GODELLE}

\textsc{Daan KRAMMER}

\textsc{Jean MICHEL}
\end{center}

\vspace{20mm}

This text consists of the introduction, table of contents, and bibliography of a long manuscript (703 pages) that is currently submitted for publication. This manuscript develops an extension of Garside's approach to braid groups and provides a unified treatment for the various algebraic structures that appear in this context.
The complete text can be found at

{\tt http://www.math.unicaen.fr/$\sim$garside/Garside.pdf}.

\noindent Comments are welcome.

\section*{Introduction}

\subsubsection*{A natural, but slowly emerging program}

In his PhD thesis prepared under the supervision of Graham Higman and defended in 1965~\cite{GarPhD}, and in the article that followed~\cite{Gar}, Frank A.\,Garside (1915--1988) solved the Conjugacy Problem of Artin's braid group~$B_n$ by introducing a submonoid~$\BP{n}$ of~$B_n$ and a distinguished element~$\Delta_n$ of~$\BP{n}$ that he called fundamental and showing that every element of~$B_n$ can be expressed as a fraction of the form~$\Delta_n^m g$ with~$m$ an integer and $g$ an element of~$\BP{n}$. Moreover, he proved that any two elements of the monoid~$\BP{n}$ admit a least common multiple, thus somehow extending to the non-Abelian groups~$B_n$ some of the standard tools available in a torsion-free Abelian group~$\ZZZZ^n$.

In the beginning of the 1970's, it was soon realized  by Brieskorn and Saito \cite{BrS} using an algebraic approach and by Deligne \cite{Dlg} using a more geometric approach that Garside's results extend to all generalized braid groups associated with finite Coxeter groups, that is, all Artin (or, better, Artin--Tits) groups of spherical type.

The next step forward was the possibility of defining, for every element of the braid monoid~$\BP{n}$ (and, more generally, of every spherical Artin--Tits monoid) a distinguished decomposition in terms of the divisors of the fundamental element~$\Delta_n$: the point is that, if $g$ is an element of~$\BP{n}$, then there exists a (unique) greatest common divisor~$g_1$ for~$g$ and~$\Delta_n$ and, moreover $g \not= 1$ implies $g_1 \not= 1$: then $g_1$ is a distinguished fragment of~$g$ (the ``head'' of~$g$) and, if we repeat the operation with the element~$g'$ that satisfies $g = g_1 g'$, we extract the head~$g_2$ of~$g'$ and, iterating, we end up with an expression $g_1 \pdots g_p$ of~$g$ in terms of divisors of~$\Delta_n$. Although Garside was very close to such a decomposition when he proved that greatest common divisors exist in~$\BP{n}$, the result does not appear in his work explicitly, and it seems that the first explicit occurrences of such distinguished decompositions, or \emph{normal forms}, goes back to the 1980's in independent work by Adjan~\cite{AdjDel}, El Rifai and Morton~\cite{ElM}, and Thurston (circulated notes~\cite{Thu}, later appearing as Chapter~IX in the book~\cite{Eps} by~Epstein {\it et al.}). The normal form was soon used to improve Garside's solution of the Conjugacy Problem~\cite{ElM} and, extended from the monoid to the group, to serve as a paradigmatic example in the then emerging theory of automatic groups of Cannon, Thurston, and others. Sometimes called the \emph{greedy normal form}---or \emph{Garside's normal form}, or \emph{Thurston's normal form}---it became a standard tool in the investigation of braids and Artin--Tits monoids and groups from a viewpoint of geometric group theory and representation, essential in particular in Krammer's algebraic proof of the linearity of braid groups~\cite{KraLin1, KraLin2}.

In the beginning of the 1990's, it was realized by one of us that some ideas from Garside's approach to braid monoids can be applied in a different context to analyze a certain ``geometry monoid''~$\LDM$ that appears in the study of the so-called left-selfdistributivity law $x(yz) = (xy)(xz)$. In particular, the criterion used by Garside to establish that the braid monoid~$\BP{n}$ is left-cancellative (that is, $gh = gh'$ implies $h = h')$ can be adapted to~$\LDM$ and a normal form reminiscent of the greedy normal form exists---with the main difference that the pieces of the normal decompositions are not the divisors of some unique element similar to Garside's fundamental braid~$\Delta_n$, but they are divisors of elements~$\Delta_t$ that depend on some object~$t$ (actually a tree) attached to the element one wishes to decompose. The approach led to results about the exotic left-selfdistributivity law~\cite{Dez} and, more unexpectedly, about braids and their orderability when it turned out that the monoid~$\LDM$ naturally projects to the (infinite) braid monoid~$\BP\infty$ \cite{Dfa, Dfb, Dgd}.

At the end of the 1990's, following a suggestion by Luis Paris, the idea arose of listing the  abstract properties of the monoid~$\BP{n}$ and of Garside's fundamental braid~$\Delta_n$ that make the algebraic theory of~$B_n$ possible. This resulted in introducing the notions of a \emph{Garside monoid} and a \emph{Garside element}~\cite{Dfx}. In a sense, this is just a sort of reverse engineering, and proving results about the existence and the properties of the normal form essentially means checking that no assumption has been forgotten in the definition. However, it soon appeared that a number of new examples were eligible, and, specially after some cleaning of the definitions was completed~\cite{Dgk}, that the new framework was really more general than the original braid framework. The main benefit was that extending the results often resulted in discovering new improved arguments no longer relying on superfluous assumptions or on specific properties. This program turned out to be rather successful and it led to many developments by a number of different authors \cite[...]{BesD, BesF, BeC, BGG1, BGG2, BGG3, CMW, ChM, CrP, FrG, Geb, GodPar1, GodPar2, Lee, LeL, McC, Pic, Sib}. Today the study of Garside monoids is still far from complete, and many questions remain open.

However, in the meanwhile, it soon appeared that, although efficient, the framework of Garside monoids as stabilized in the 1990's is far from optimal. Essentially, several assumptions, in particular Noetherianity conditions, are superfluous and they just discard further natural examples. Also, excluding nontrivial invertible elements appears as an artificial limiting assumption. More importantly, one of us (DK) in a 2005 preprint subsequently published as~\cite{KraCel} and two of us (FD, JM)~\cite{DiMCat}, as well as David Bessis in an independent research~\cite{Bes}, realized that normal forms similar to those involved in Garside monoids can be developed and usefully applied in a context of categories, leading to what they naturally called \emph{Garside categories}. By the way, similar structures are already implicit in the 1976 paper~\cite{DeL} by Deligne--Lusztig, as well as in the above mentioned example of~$\LDM$~\cite{Dfb, Dgd}, and in EG's PhD thesis~\cite{GodPhD}.

It was therefore time around 2007 for the development of a new, unifying framework that would include all the previously defined notions, remove all unneeded assumptions, and allow for optimized arguments. This program was developed in particular during a series of workshops and meetings between 2007 and 2012, and it resulted in the current text. As suggested in the above account, the emphasis is put on the normal form and its mechanism, and the framework is that of a general category with only one assumption, namely left-cancellativity. Then the central notion is that of a \emph{Garside family}, defined to be any family that gives rise to a normal form of the expected type. Then, of course, every Garside element~$\Delta$ in a Garside monoid provides an example of a Garside family, namely the set of all divisors of~$\Delta$, but many more Garside families may exist---and they do, as we shall see in the text. Note as, in a sense, our current generalization is the ultimate one since, by definition, no further extension may preserve the existence of a greedy normal form.  However, different approaches might be developed, either by relaxing the definition of a greedy decomposition (see the Notes at the end of Chapter~III) or, more radically, by putting the emphasis on other aspects of Garside groups rather than on normal forms. Typically, several authors, including J.\,Crisp, J.\,McCammond and one of us (DK) proposed to view a Garside group mainly as a group acting on a lattice in which certain intervals of the form~$[1, \Delta]$ play a distinguished role, thus paving the way for other types of extensions.

Our hope---and our claim---is that the new framework so constructed is quite satisfactory. By this, we mean that most of the properties previously established in more particular contexts can be extended to larger contexts. It is \emph{not} true that all properties of, say, Garside monoids extend to arbitrary categories equipped with a Garside family but, in most cases, addressing the question in an extended framework helps improving the arguments and really capturing the essential features. Typically, almost all known properties of Garside monoids do extend to categories that admit what we call a bounded Garside family, and the proofs cover for free all previously considered notions of Garside categories.

It is clear that a number future developments will continue to involve particular types of monoids or categories only: we do not claim that our approach is universal... However, we would be happy if the new framework---and the associated terminology---could become a natural reference for further works.

\subsubsection*{About this text}

 The aim of the current text is to give a state-of-the-art presentation of this approach. Finding a proper name turned out to be not so obvious. On the one hand, ``Garside calculus'' would be a natural title, as the greedy normal form and its variations are central in this text: although algorithmic questions are not emphasized, most constructions are effective and the mechanism of the normal form is indeed a sort of calculus. On the other hand however, many results, in particular those of structural nature, exploit the normal form but are not reducible to it, making a title like ``Garside structures'' or ``Garside theory" more appropriate. But such a title is certainly too ambitious for what we can offer: no genuine structure theory or no exhaustive classification of, say, Garside families is to be expected at the moment. What we do here is to develop a framework that, we think and hope, can become a good base for a still-to-come theory.  Another option could have been ``Garside categories'', but it will be soon observed that no notion with that name is introduced here: in view of the subsequent developments, a reasonable meaning could be ``a cancellative category that admits a Garside map'', but a number of variations are still possible, and any particular choice could become obsolete soon---as is, in some sense, the notion of a Garside group.  So, finally, our current title, ``Foundations of Garside Theory'', may be the one that reflects the current content in the best way: the current text should be seen as an invitation for further research, and does not aim at being exhaustive---reporting about all previous results involving Garside structures would already be very difficult---but concentrates on what seems to be the core of the subject. 
 
There are two parts. Part~A is devoted to general results, and it offers a very careful treatment of the bases. Here complete proofs are given, and the results are illustrated with a few basic examples. By contrast, Part~B consists of essentially independent chapters explaining further examples or families of examples that are in general more elaborate. Here some proofs may be omitted, and the discussion is centered around what can be called the Garside aspects in the considered structures. 

 Our general scheme will be to start from an analysis of normal decompositions and then to introduce Garside families as the framework guaranteeing the existence of normal decompositions. Then the three main questions we shall address and a chart of the corresponding chapters looks as follows:
$$\begin{tabular}{l}
\hline
\VR(4,2)$\bullet$ {\bf How} do Garside structures work? (mechanism of normal decomposition)\\
\HS{8}Chapter~III (domino rules, geometric aspects)\\
\HS{8}Chapter~VII (compatibility with subcategories)\\
\VR(0,3)\HS{8}Chapter~VIII (connection with conjugacy)\\
\hline
\VR(4,2)$\bullet$ {\bf When} do Garside structures exist? (existence of normal decomposition)\\
\HS{8}Chapter~IV (recognizing Garside families)\\
\HS{8}Chapter~VI (recognizing Garside germs)\\
\VR(0,3)\HS{8}Chapter~V (recognizing Garside maps)\\
\hline
\VR(4,2)$\bullet$ {\bf Why} consider Garside structures? (examples and applications)\\
\HS{8}Chapter~I (basic examples)\\
\HS{8}Chapter~IX (braid groups)\\
\HS{8}Chapter~X (Deligne--Luzstig varieties)\\
\HS{8}Chapter~XI (selfdistributivity)\\
\HS{8}Chapter~XII (ordered groups)\\
\HS{8}Chapter~XIII (Yang--Baxter equation)\\
\VR(0,3)\HS{8}Chapter~XIV (four more examples)\\
\hline
\end{tabular}$$
Above, and in various places, we use ``Garside structure'' as a generic and informal way to refer to the various objects occurring with the name ``Garside'': Garside families, Garside groups, Garside maps, \emph{etc}. 

\subsubsection*{The chapters}

To make further reference easy, each chapter in Part~A begins with a summary of the main results. At the end of each chapter, exercises are proposed, and a note section provides historical references, comments, and questions for further research.

Chapter~I is introductory and lists a few examples. The chapter starts with some classical examples of Garside monoids, such as free Abelian monoids or classical and dual braid monoids, and it continues with some examples of structures that are not Garside monoids but nevertheless possess a normal form similar to that of Garside monoids, thus providing a motivation for the construction of a new, extended framework.

Chapter~II is another introductory chapter in which we fix some terminology and basic results about categories and derived notions, in particular connected with divisibility relations that play an important r\^ole in the sequel. A few general results about Noetherian categories and groupoids of fractions are established. The final section describes an general method called reversing for investigating a presented category. As the question is not central in our current approach (and although it owes much to Garside's methods), some proofs of this section are deferred to an appendix at the end of the book.

Chapter~III is the one where the theory really starts. Here the notion of a normal decomposition is introduced, as well as the notion of a Garside family, abstractly introduced as a family that guarantees the existence of an associated normal form. The mechanism of the normal form is analyzed, both in the case of a category (``positive case'') and in the case of its enveloping groupoid (``signed case''): some simple diagrammatic patterns, the domino rules, are crucial, and their local character directly implies various geometric consequences, in particular a form of automaticity and the Grid Property, a strong convexity statement.

Chapter~IV is devoted to obtaining concrete characterizations of Garside families, hence, in other words, to describing assumptions that guarantee the existence of normal decompositions. In this chapter, one establishes external characterizations, meaning that we start with a category~$\CCC$ and look for conditions ensuring that a given subfamily~$\SSS$ of~$\CCC$ is a Garside family. Various answers are given, in a general context first, and then in particular contexts where some conditions come for free: typically, if the ambient category~$\CCC$ is Noetherian and admits unique least common right-multiples, then a subfamily~$\SSS$ of~$\CCC$ is a Garside family if and only if it generates~$\CCC$ is is closed under least common right-multiple and right-divisor.

Chapter~V investigates particular Garside families that are called bounded. Essentially, a Garside family~$\SSS$ is bounded is there exists a map~$\Delta$ (an element in the case of a monoid) such that $\SSS$ consists of the divisors of~$\Delta$ (in some convenient sense). Not all Garside families are bounded, and, contrary to the existence of a Garside family, the existence of a bounded Garside family is not guaranteed in every category. Here we show that a bounded Garside family is sufficient to prove most of the results previously established for a Garside monoid, including the construction of $\Delta$-normal decompositions, a variant of the symmetric normal decompositions used in groupoids of fractions.

Chapter~VI provides what can be called internal (or intrinsic) characterizations of Garside families: here we start with a family~$\SSS$ equipped with a partial product, and we wonder whether there exists a category~$\CCC$ in which $\SSS$ embeds as a Garside family. The good news is that such characterizations do exist, meaning that, when the conditions are satisfied, all properties of the generated category can be read inside the initial family~$\SSS$. This local approach turns to be very useful to construct examples and, in particular, it can be used to construct a sort of unfolded, torsion-free version of convenient groups, typically braid groups starting from Coxeter groups.

Chapter~VII is devoted to subcategories. Here one investigates natural questions such as the following: if $\SSS$ is a Garside family in a category~$\CCC$ and $\CCCu$ is a subcategory of~$\CCC$, then is $\SSS \cap \CCCu$ a Garside family in~$\CCCu$ and, if so, what is the connection between the associated normal decompositions? Of particular interest are the results involving subgerms, which somehow provide a possibility of reading inside a given Garside family~$\SSS$ the potential properties of the subcategories generated by the subfamilies of~$\SSS$.

Chapter~VIII addresses conjugacy, first in the case of a category equipped with an arbitrary Garside family, and then, mainly, in the case of a category equipped with a bounded Garside family. Here again, most of the results previously established for Garside monoids can be extended, including the cycling, decycling, and sliding transformations which provide a decidability result for the Conjugacy Problem whenever convenient finiteness assumptions are satisfied. We also extend the geometric methods of Bestvina to describe periodic elements in this context.

Part~B begins with Chapter~IX devoted to (generalized) braid groups. Here we show how both the reversing approach of Chapter~II and the germ approach of Chapter~VI can be applied to construct and analyze classical and dual Artin--Tits monoids. We also briefly mention the braid groups associated with complex reflection groups, as well as several exotic Garside structures on~$B_n$. The applications of Garside structures in the context of braid groups are too many to be described exhaustively, and we just list some of them in the Notes section.

Chapter~X is a direct continuation of Chapter~IX. It reports about the use of Garside-type methods in the study of Deligne--Lusztig varieties, an ongoing program that aims at establishing by a direct proof some of the consequences of the Brou\'e Conjectures about finite reductive groups. Several questions in this approach directly involve conjugacy in generalized braid groups, and the results of Chapter~VIII are then crucial.

Chapter~XI is an introduction to the Garside structure hidden in the above mentioned algebraic law $x(yz) = (xy)(xz)$, a typical example where a categorical framework is needed (or, at the least, the framework of Garside monoids is not sufficient). Here a promising contribution of the Garside approach is a natural program possibly leading to the so-called Embedding Conjecture, a deep structural result that resisted all attempts so far.

Chapter~XII develops an approach to ordered groups based on divisibility properties and Garside elements, resulting in the construction of groups with the property that the associated space of orderings contains isolated points, which answers one of the natural questions of the area. Braid groups are typical examples, but considering what we call triangular presentations leads to a number of different examples.  

Chapter~XIII is a self-contained introduction to set-theoretic solutions of the Yang--Baxter equation and the associated structure groups, which make an important family of Garside groups. The exposition is centered on the connection between the RC-law $(xy)(xz) = (yx)(yz)$ and the right-complement operation on the one hand, and what is called the geometric $I$-structure on the other hand. Here the Garside approach both provides a specially efficient framework, in particular for reproving results about the RC-law, and leads to new results.

Chapter~XIV presents four  unrelated topics involving interesting Garside families: divided categories and decompositions categories with two applications, then an extension of the framework of Chapter~XIII to more general RC-systems, then what is called the braid group of~$\ZZZZ^n$, a sort of analog of Artin's braid group in which permutations of $\{1 \wdots n\}$ are replaced with linear orderings of~$\ZZZZ^n$, and, finally, an introduction to groupoids of cell decompositions that arise when the mapping class group approach to braid groups is extended by introducing sort of roots of the generators~$\sig{i}$. 

The final Appendix contains the postponed proofs of some technical statements from Chapter~II for which no complete reference exists in literature.

\subsubsection*{Thanks}

We primarily wish to thank David Bessis, who was part of the crew at an early stage of the project, but then quitted it for personal reasons.

Next, we thank all the colleagues and students who participated in the various meetings dedicated to the project and who contributed by their questions and suggestions to this text. A (certainly non-exhaustive) list includes Marc Autord, Serge Bouc, Michel Brou\'e, Matthieu Calvez, Ruth Corran, Jean Fromentin, Volker Gebhardt, Tomas Gobet, Juan Gonz\'alez-Meneses, Tatiana Ivanova--Gateva, Jean-Yves H\'ee, Eric Jespers, Eon-Kyung Lee, Sang-Jin Lee, Jon McCammond, Ivan Marin, Jan Okni\'nski, Luis Paris, Matthieu Picantin, Maya Van Campenhout, Bert Wiest. Also, special thoughts for Joan Birman and Hugh Morton, whose support and interest in the subject has always been strong.

Let us mention that the project was supported partly by the ANR grant TheoGar ANR-08-BLAN-0269-02.

\bigskip
\begin{flushright}
Caen, Amiens, Warwick, Paris, June 2013

\medskip

Patrick Dehornoy

Fran\c cois Digne

Eddy Godelle

Daan Krammer

Jean Michel
\end{flushright}

\vfill

\begin{center}
\emph{What remains to be done?}
\end{center}

\bigskip

$\bullet$ Add a few more examples in the text, typically involving the wreathed free Abelian monoid~$\widetilde\NNNN_n$ in addition to the standard ones involving the braid group~$B_n$ in order to show how the results and algorithms look like when there are nontrivial invertible elements.

$\bullet$ Uniformize the visual aspect of all pictures (same style of arrows, same linewidth, etc.);

$\bullet$ Post the solutions of the exercises (which are written but not printed here) on a dedicated website;

$\bullet$ Maybe : transform some secondary statements that are not subsequently referred to into exercises.

\newpage

\begin{center}
\Large{\sc Contents}
\end{center}

\Chapter{Introduction}vii

\Part{A. General theory}1

\Chapter{I. Some examples}3

\Section{1. Classical examples}3

\Subsection{1.1. Free Abelian monoids}3

\Subsection{1.2. Braid groups and monoids}5

\Subsection{1.2. Dual braid monoids}9

\Section{2. Garside monoids and groups}11

\Subsection{2.1. The notion of a Garside monoid}11

\Subsection{2.2. More examples}13

\Section{3. Why a further extension?}15

\Subsection{3.1. Infinite braids}15

\Subsection{3.2. The Klein bottle group}16

\Subsection{3.3. Wreathed free Abelian groups}18

\Subsection{3.4. Ribbon categories}19

\Exercises22

\Notes22

\Chapter{II. Preliminaries}27

\Section{1. The category context}29

\Subsection{1.1. Categories and monoids}29

\Subsection{1.2. Subfamilies and subcategories}31

\Subsection{1.3. Invertible elements}33

\Subsection{1.4. Presentations}37

\Section{2. Divisibility and Noetherianity}40

\Subsection{2.1. Divisibility relations}41

\Subsection{2.2. Lcms and gcds}42

\Subsection{2.3. Noetherianity conditions}46

\Subsection{2.4. Height}53

\Subsection{2.5. Atoms}57

\Section{3. Groupoids of fractions}60

\Subsection{3.1. The enveloping groupoid of a category}60

\Subsection{3.2. Groupoid of fractions}62

\Subsection{3.3. Ore subcategories}65

\Subsection{3.4. Torsion elements in a groupoid of fractions}66

\Section{4. Working with presented categories}67

\Subsection{4.1. A toolbox}67

\Subsection{4.2. Right-reversing: definition}71

\Subsection{4.3. Right-reversing: termination}75

\Subsection{4.4. Right-reversing: completeness}79

\Exercises87

\Notes90

\Chapter{III. Normal decompositions}93

\Section{1. Greedy decompositions}95

\Subsection{1.1. The notion of an $\mathcal{S}$-greedy path}96

\Subsection{1.2. The notion of an $\mathcal{S}$-normal path}100

\Subsection{1.3. The notion of a Garside family}105

\Subsection{1.4. Recognizing Garside families}107

\Subsection{1.5. The second domino rule}114

\Section{2. Symmetric normal decompositions}118

\Subsection{2.1. Left-disjoint elements}118

\Subsection{2.2. Symmetric normal decompositions}121

\Subsection{2.3. Uniqueness of symmetric normal decompositions}126

\Subsection{2.4. Existence of symmetric normal decompositions}123

\Subsection{2.5. More domino rules}132

\Subsections{Appendix: existence of symmetric normal decompositions,}

\Subsection{\indent\indent general case of a left-cancellative category}136

\Section{3. Geometric and algorithmic properties}139

\Subsection{3.1. Geodesics}140

\Subsection{3.2. The Grid Property}141

\Subsection{3.3. The Fellow Traveller Property}146

\Subsection{3.4. The Garside resolution}152

\Subsection{3.5. Word Problem}161

\Exercises164

\Notes165

\Chapter{IV. Garside families}169

\Section{1. The general case}172

\Subsection{1.1. Closure properties}172

\Subsection{1.2. Characterizations of Garside families}180

\Subsection{1.3. Special Garside families}183

\Subsection{1.4. Head functions}187

\Section{2. Special contexts}191

\Subsection{2.1. Solid families}191

\Subsection{2.2. Right-Noetherian categories}194

\Subsection{2.3. Categories that admit right-mcms}200

\Subsection{2.4. Categories with unique right-lcms}207

\Subsection{2.5. Finite height}209

\Section{3. Geometric and algorithmic applications}211

\Subsection{3.1. Presentations}211

\Subsection{3.2. Isoperimetric inequalities}215

\Subsection{3.3. Word Problem}216

\Subsection{3.4. Case of categories that admit lcms}219

\Exercises225

\Notes227

\Chapter{V. Bounded Garside families}231

\Section{1. Right-bounded Garside families}234

\Subsection{1.1. The notion of a right-bounded Garside family}234

\Subsection{1.2. Right-Garside maps}237

\Subsection{1.3. The Garside functor $\phi_\Delta$}240

\Subsection{1.4. Powers of a right-bounded Garside family}244

\Subsection{1.5. Preservation of normality}247

\Section{2. Bounded Garside families}250

\Subsection{2.1. The notion of a bounded Garside family}251

\Subsection{2.2. Powers of a bounded Garside family}254

\Subsection{2.3. The case of a cancellative category}256

\Subsection{2.4. Garside maps}259

\Subsection{2.5. Existence of lcms and gcds}262

\Section{3. Delta-normal decompositions}266

\Subsection{3.1. The positive case}266

\Subsection{3.2. The general case}270

\Subsection{3.3. Symmetric normal decompositions}277

\Subsection{3.4. Co-normal decompositions}279

\Exercises282

\Notes283

\Chapter{VI. Germs}287

\Section{1. Germs}290

\Subsection{1.1. The notion of a germ}290

\Subsection{1.2. The embedding problem}292

\Subsection{1.3. Atoms in a germ}295

\Subsection{1.4. Garside germs}297

\Section{2. Recognizing Garside germs}300

\Subsection{2.1. The families $\mathscr{I}_{\mathcal{S}}$ and $\mathscr{J}_{\mathcal{S}}$}300

\Subsection{2.2. Greatest $\mathscr{I}$-functions}307

\Subsection{2.3. Noetherian germs}310

\Subsection{2.4. An application: germs derived from a groupoid}313

\Section{3. Bounded germs}319

\Subsection{3.1. Right-bounded germs}319

\Subsection{3.2. Bounded germs}321

\Subsection{3.3. An application: germs from lattices}324

\Exercises326

\Notes327

\Chapter{VII. Subcategories}329

\Section{1. Subcategories}332

\Subsection{1.1. Closure under quotient}332

\Subsection{1.2. Subcategories that are closed under $=^{\!\scriptscriptstyle\times}$}336

\Subsection{1.3. Head-subcategories}338

\Subsection{1.4. Parabolic subcategories}344

\Section{2. Compatibility with a Garside family}345

\Subsection{2.1. Greedy paths}346

\Subsection{2.2. Compatibility with a Garside family}347

\Subsection{2.3. Compatibility, special subcategories}351

\Subsection{2.4. Compatibility with symmetric decompositions}354

\Section{3. Subfamilies of a Garside family}357

\Subsection{3.1. Subgerms}358

\Subsection{3.2. Transfer results}361

\Subsection{3.3. Garside subgerms}365

\Subsection{3.4. Head-subgerms}371

\Section{4. Subcategories associated with functors}372

\Subsection{4.1. Subcategories of fixed points}372

\Subsection{4.2. Image subcategory}374

\Exercises379

\Notes381

\Chapter{VIII. Conjugacy}385

\Section{1. Conjugacy categories}387

\Subsection{1.1. General conjugacy}387

\Subsection{1.2. Cyclic conjugacy}392

\Subsection{1.3. Twisted conjugacy}396

\Subsection{1.4. An example: ribbon categories}400

\Section{2. Cycling, sliding, summit sets}409

\Subsection{2.1. Cycling and decycling}409

\Subsection{2.2. Sliding circuits}417

\Section{3. Conjugacy classes of periodic elements}429

\Subsection{3.1. Periodic elements}429

\Subsection{3.2. Geometric methods}431

\Subsection{3.3. Conjugates of periodic elements}438

\Exercises443

\Notes443

\Part{B. Specific examples}447

\Chapter{IX. Braids}449

\Section{1. The classical Garside structure on Artin--Tits groups}449

\Subsection{1.1. Coxeter groups}450

\Subsection{1.2. Artin--Tits groups, reversing approach}455

\Subsection{1.3. Artin--Tits groups, germ approach}461

\Section{2. More Garside structures on Artin--Tits groups}463

\Subsection{2.1. The dual braid monoid}463

\Subsection{2.2. The case of the symmetric group}465

\Subsection{2.3. The case of finite Coxeter groups}469

\Subsection{2.4. Exotic Garside structures on $B_n$}470

\Section{3. Braid groups of well-generated complex reflection groups}473

\Subsection{3.1. Complex reflection groups}473

\Subsection{3.2. Braid groups of complex reflection groups}475

\Subsection{3.3. Well-generated complex reflection groups}476

\Subsection{3.4. Tunnels}477

\Subsections{3.5. The Lyashko--Looijenga covering and Hurwitz action}

\Subsection{\indent\indent on decompositions of~$\delta$}479

\Exercises482

\Notes482

\Chapter{X. Deligne--Lusztig varieties}487

\Section{1. Finite linear groups as reductive groups}487

\Subsection{1.1. Reductive groups}487

\Subsection{1.2. Some important subgroups}488

\Subsection{1.3. $\boldsymbol{G}^F$-conjugacy}489

\Section{2. Representations}491

\Subsection{2.1. Complex representations of $\boldsymbol{G}^F$}491

\Subsection{2.2. Deligne--Lusztig varieties}492

\Subsection{2.3. Modular representation theory}493

\Section{3. Geometric Brou\'e Conjecture, torus case}494

\Subsection{3.1. The geometric approach}495

\Subsection{3.2. Endomorphisms of Deligne--Lusztig varieties}496

\Section{4. Geometric Brou\'e Conjecture, the general case}499

\Subsection{3.1. The parabolic case}500

\Subsection{3.2. The really general case}503

\Notes506

\Chapter{XI. Left self-distributivity}509

\Section{1. Garside sequences}510

\Subsection{1.1. Partial actions}510

\Subsection{1.2. Right-Garside sequences}513

\Subsection{1.3. Derived notions}515

\Section{2. LD-expansions and the category $\mathcal{L\!D}_0$}518

\Subsection{2.1. Free LD-systems}518

\Subsection{2.2. LD-expansions}520

\Subsection{2.3. The category $\mathcal{L\!D}_0$}522

\Subsection{2.4. Simple LD-expansions}524

\Section{3. Labelled LD-expansions and the category $\mathcal{L\!D}$}524

\Subsection{3.1. The operators $\Sigma_\alpha$}525

\Subsection{3.2. The monoid $M_{\scriptscriptstyle L\!D}$}526

\Subsection{3.3. The category $\mathcal{L\!D}$}529

\Subsection{3.4. The Embedding Conjecture}532

\Section{4. Connection with braids}536

\Subsection{4.1. The main projection}536

\Subsection{4.2. Reproving braid properties}538

\Subsection{4.2. Hurwitz action of braids on LD-systems}542

\Exercises544

\Notes545

\Chapter{XII. Ordered groups}549

\Section{1. Ordered groups and monoids of $O$-type}549

\Subsection{1.1. Orderable and bi-orderable groups}550

\Subsection{1.2. The space of orderings on a group}553

\Subsection{1.3. Two examples}556

\Section{2. Construction of isolated orderings}558

\Subsection{2.1. Triangular presentations}558

\Subsection{2.2. Existence of common multiples}561

\Subsection{2.3. More examples}564

\Subsection{2.4. Effectivity questions}566

\Section{3. Further results}568

\Subsection{3.1. Dominating elements}569

\Subsection{3.2. Right-ceiling}570

\Subsection{3.3. The specific case of braids}571

\Exercises574

\Notes575

\Chapter{XIII. Set-theoretic solutions of YBE}579

\Section{1. Several equivalent frameworks}579

\Subsection{1.1. Set-theoretic solutions of the Yang--Baxter equation}580

\Subsection{1.2. Involutive biracks}582

\Subsection{1.3. RC- and RLC-quasigroups}584

\Section{2. Structure monoids and groups}589

\Subsection{2.1. Structure monoids and groups}590

\Subsection{2.2. RC-calculus}593

\Subsection{2.3. Every structure monoid is a Garside monoid}598

\Subsection{2.4. A converse connection}600

\Section{3. Monoids of $I$-type}603

\Subsection{3.1. The $I$-structure}603

\Subsection{3.2. Monoids of $I$-type}606

\Subsection{3.3. Coxeter-like groups}610

\Exercises615

\Notes615

\Chapter{XIV. More examples}619

\Section{1. Divided and decomposition categories}619

\Subsection{1.1. Divided categories}620

\Subsection{1.2. Decomposition categories}625

\Section{2. Cyclic systems}632

\Subsection{2.1. Weak RC-systems}632

\Subsection{2.2. Units and ideals}635

\Subsection{2.3. The structure category of a weak RC-system}638

\Section{3. The braid group of $\mathbb{Z}^n$}644

\Subsection{3.1. Ordering orders}644

\Subsection{3.2. Lexicographic orders of $\mathbb{Z}^n$}645

\Subsection{3.3. A lattice ordering on $\mathrm{GL}(n,\mathbb{Z})$}647 

\Section{4. Cell decompositions of a punctured disk}649

\Subsection{4.1. Braid groups as mapping class groups}649

\Subsection{4.2. Cell decompositions}651

\Subsection{4.3. The group $B_\ell$ and the category $\mathcal{B}_\ell$}652 

\Subsection{4.4. Flips}654

\Subsection{4.5. A bounded Garside family}658 

\Exercises659

\Notes660

\Chapter{Appendix: Some missing proofs for Chapter II}663

\Section{1. Groupoid of fractions}663

\Subsection{1.1. Ore's theorem}663

\Subsection{1.2. Ore subcategories}668

\Section{2. Working with presented categories}669

\Subsection{2.1. Right-reversing: termination}669

\Subsection{2.2. Right-reversing: completeness}670

\Exercises676

\Chapter{Bibliography}677

\Chapter{Index}685

\newpage

\newcommand{\noopsort}[1]{}
\providecommand{\bysame}{\leavevmode\hbox to3em{\hrulefill}\thinspace}
\providecommand{\MR}{\relax\ifhmode\unskip\space\fi MR }
\providecommand{\MRhref}[2]{%
  \href{http://www.ams.org/mathscinet-getitem?mr=#1}{#2}
}
\providecommand{\href}[2]{#2}

\vfill

Addresses of the authors:\smaller

\bigskip\sc PD: {Laboratoire de Math\'ematiques Nicolas Oresme,
CNRS UMR 6139, Universit\'e de Caen, 14032 Caen, France}

{\tt patrick.dehornoy@unicaen.fr \quad www.math.unicaen.fr/\!\hbox{$\sim$}dehornoy}

\bigskip\sc FD: {Laboratoire Ami\'enois de Math\'ematique Fondamentale et Appliqu\'ee, 
CNRS UMR 7352, Universit\'e de Picardie
Jules-Verne, 80039 Amiens, France}

{\tt digne@u-picardie.fr \quad www.mathinfo.u-picardie.fr/digne/}

\bigskip\sc EG: {Laboratoire de Math\'ematiques Nicolas Oresme,
CNRS UMR 6139, Universit\'e de Caen, 14032 Caen, France}

{\tt godelle@math.unicaen.fr \quad www.math.unicaen.fr/\!\hbox{$\sim$}godelle}

\bigskip\sc DK: {Mathematics Institute, University of Warwick, Coventry CV4
7AL, United Kingdom}

{\tt D.Krammer@warwick.ac.uk \quad www.warwick.ac.uk/\!\hbox{$\sim$}masbal/}

\bigskip\sc JM: {Institut de Math\'ematiques de Jussieu, CNRS UMR 7586, 
Universit\'e Denis Diderot Paris 7, 75205 Paris 13, France}

{\tt jmichel@math.jussieu.fr \quad www.math.jussieu.fr/\!\hbox{$\sim$}jmichel/}
\end{document}